\documentclass[11pt]{amsart}
\usepackage[a4paper, margin=2cm, headsep=0.7cm, footskip=1cm]{geometry}             
\usepackage[dvipsnames]{xcolor}   		             		
\usepackage{graphicx}				
\usepackage{amssymb,mathrsfs}
\usepackage{amsthm,amsmath,stmaryrd}
\usepackage{tikz}
\usepackage{tikz-cd}
\usepackage{accents,upgreek, enumitem}
\usepackage{bm}
\usepackage{mathtools}
\usepackage[all]{xy}
\usepackage{caption}
\usepackage{array,booktabs}

\tikzset{
  commutative diagrams/.cd, 
  arrow style=tikz, 
  diagrams={>=stealth}
}

\DeclareMathAlphabet{\mathpzc}{OT1}{pzc}{m}{it}

\usepackage{hyperref}
\hypersetup{
  colorlinks   = true,          
  urlcolor     = blue,          
  linkcolor    = purple,          
  citecolor   = blue             
}

\usepackage[capitalize]{cleveref}

\theoremstyle{plain}
\newtheorem*{theorem}{Theorem}
\newtheorem*{example}{Example}

\theoremstyle{definition}
\newtheorem*{remark}{Remark}
\newtheorem*{acknowledgements}{Acknowledgements}

\title[Hamiltonian isotopies are orientation-preserving]{Hamiltonian isotopies of relatively exact Lagrangians \\ are orientation-preserving}

\author{Jack Smith}
\address{St John's College, Cambridge, CB2 1TP, United Kingdom}
\email{\href{mailto:j.smith@dpmms.cam.ac.uk}{j.smith@dpmms.cam.ac.uk}}

\newcommand{\HLM}[1]{\mathcal{H}_L^{\smash{#1}}}
\renewcommand{\phi}{\varphi}
\newcommand{\rH}{\operatorname{H}}
\newcommand{\HF}{\operatorname{HF}}

\newcommand{\rmd}{\mathrm{d}}
\newcommand{\Ham}{\operatorname{Ham}}
\newcommand{\Symp}{\operatorname{Symp}}
\newcommand{\id}{\operatorname{id}}
\newcommand{\GL}{\operatorname{GL}}

\newcommand{\CC}{\mathbb{C}}
\newcommand{\PP}{\mathbb{P}}
\newcommand{\RR}{\mathbb{R}}
\newcommand{\ZZ}{\mathbb{Z}}

\setenumerate{label=(\alph*)}

\begin{document}

\begin{abstract}
Given a closed, orientable Lagrangian submanifold $L$ in a symplectic manifold $(X, \omega)$, we show that if $L$ is relatively exact then any Hamiltonian diffeomorphism preserving $L$ setwise must preserve its orientation.  In contrast to previous results in this direction, there are no spin hypotheses on $L$.  Curiously, the proof uses only mod-$2$ coefficients in its singular and Floer cohomology rings.
\end{abstract}

\maketitle

Let $(X, \omega)$ be a symplectic manifold which is compact or well-behaved (e.g.~convex and cylindrical) at infinity.  Recall that an isotopy $\phi_t$ of $X$ is \emph{Hamiltonian} if is generated by a time-dependent smooth function $H_t : X \to \RR$, in the sense that the generating vector field $V_t = \phi_t^* \dot{\phi}_t$ satisfies
\[
\omega(V_t, {-}) = - \rmd H_t.
\]
The group $\Ham(X)$ comprises the end-points $\phi_1$ of Hamiltonian isotopies $\phi_t$ starting from $\phi_0 = \id_X$, and is an important subgroup of the symplectomorphism group $\Symp(X)$.

Let $L \subset X$ be a closed, orientable Lagrangian submanifold, and define
\[
\Ham(X, L) = \{\phi \in \Ham(X) : \phi(L) = L\}.
\]
Assume throughout that $L$ is relatively exact, meaning that $\omega$ vanishes on $\pi_2(X, L)$.

The purpose of this short note is to prove the following result.

\begin{theorem}
The action of $\Ham(X, L)$ on $L$ is orientation-preserving.
\end{theorem}

\begin{example}
\label{ex}
If $L$ is homotopy equivalent to $\RR\PP^{2n+1}$ then $\Ham(X, L)$ acts trivially on $\rH_*(L; \ZZ)$.
\end{example}

\begin{remark}
The $x$-axis in $\RR^2$, with its standard symplectic form, can be rotated through angle $\pi$ by a Hamiltonian isotopy, so compactness of $L$ is essential in the theorem.  Similarly, the equator in $S^2$ can be turned over by a Hamiltonian isotopy, so relative exactness is also essential.
\end{remark}

The action of $\Ham(X, L)$ on $L$ induces a homomorphism
\[
\Ham(X, L)  \to \GL(\rH_*(L; \ZZ)),
\]
whose image is the \emph{(Hamiltonian) homological Lagrangian monodromy group} $\HLM{\ZZ}$ of $L$.  The theorem can be interpreted as saying that $\HLM{\ZZ}$ fixes the fundamental class $[L]$, whilst in the example $\HLM{\ZZ}$ is trivial.

The group $\HLM{\ZZ}$ was introduced by M.-L.~Yau in \cite{YauMonodromyAndIsotopy} to study the monotone Clifford and Chekanov tori in $\CC^2$.  It was extensively studied in the weakly exact setting, along with its analogue $\HLM{\ZZ/2}$ in $\GL(\rH_*(L; \ZZ/2))$, by Hu--Lalonde--Leclercq in \cite{HuLalondeLeclercq}.  They showed that $\HLM{\ZZ/2}$ is trivial in this case, and that $\HLM{\ZZ}$ is also trivial if $L$ is (relatively) spin and one restricts to those $\phi \in \Ham(X, L)$ preserving the (relative) spin structure.  Porcelli \cite{Porcelli} subsequently gave a new proof of these results under weaker hypotheses, which extends to generalised homology theories.  In particular, he showed that $\HLM{\ZZ}$ is trivial if $L$ is spin, with no restriction on $\phi$.  Compared with these works, our result is therefore notable because of its lack of spin hypotheses on both $L$ and $\phi$.  Note that the $\RR\PP^{2n+1}$ example is not spin for any even $n \geq 2$.

\begin{proof}[Proof of the theorem]
We work throughout with $\ZZ/2$ coefficients.  First we recall some Floer theory.  Suppose $L_1$ and $L_2$ are closed, orientable Lagrangian submanifolds in $X$ for which the Floer cohomology $\HF(L_1, L_2)$ can be defined.  This holds, for example, if they are relatively exact and Hamiltonian isotopic to each other.  For each $L_i$ a choice of orientation corresponds to a $\ZZ/2$-grading in the sense of \cite{SeidelGraded}, with respect to the canonical $\ZZ/2$-grading of $X$.  Orienting both $L_1$ and $L_2$ then gives a $\ZZ/2$-grading on $\HF^*(L_1, L_2)$.  Reversing the orientation of $L_1$ corresponds to shifting its grading, which we denote by $L_1[1]$, and there is then a natural identification $\HF^*(L_1[1], L_2) = \HF^{*-1}(L_1, L_2)$.  Similarly for $L_2$, with $\HF^*(L_1, L_2[1]) = \HF^{*+1}(L_1, L_2) = \HF^{*-1}(L_1, L_2)$.

Given a third closed, oriented Lagrangian $L_3$ for which the Floer cohomology with $L_1$ and $L_2$ can be defined, there is a Floer product
\[
\mu^2 : \HF^*(L_2, L_3) \otimes \HF^*(L_1, L_2) \to \HF^*(L_1, L_3).
\]
This respects the $\ZZ/2$-gradings, and is associative in the obvious sense if we also consider a fourth Lagrangian.  Each $\HF^*(L_i, L_i)$ then becomes a $\ZZ/2$-graded ring, and this ring is unital.  It is well-known that if $L_i$ is relatively exact then the ring $\HF^*(L_i, L_i)$ is given by $\rH^*(L_i)$ with the grading collapsed modulo $2$.

Now focus on our single closed, relatively exact Lagrangian $L$, on which we fix an orientation.  A compactly-supported Hamiltonian isotopy $\phi_t$ induces a \emph{continuation element}
\[
c(\phi_t) \in \HF^0(\phi_0(L), \phi_1(L)).
\]
This is functorial, in the sense that Floer products of continuation elements correspond to concatenations of Hamiltonian isotopies, and if $\phi_t$ is a constant isotopy then $c(\phi_t)$ is the unit.  In particular this means that continuation elements are invertible with respect to the Floer product.

Suppose for contradiction that there exists a $\phi \in \Ham(X, L)$ which reverses the orientation of $L$, and fix a Hamiltonian isotopy $\phi_t$ from $\id_X$ to $\phi$.  By cutting off the generating function $H_t$ outside a compact set containing the sweepout $\bigcup_t \phi_t(L)$ of $L$, we may assume that the isotopy is compactly supported.  It therefore has an associated continuation element $c(\phi_t) \in \HF^0(L, \phi(L))$.

Since $\phi(L)$ is $L$ with its orientation reversed, i.e.~$\phi(L) = L[1]$, we can view $c(\phi_t)$ as an element of $\HF^1(L, L)$.  By the above discussion, this element is invertible, and $\HF^*(L, L)$ is $\rH^*(L)$ with grading collapsed mod $2$, so we conclude that there exists an invertible element of odd degree in $\rH^*(L)$.  But this is impossible: any odd-degree element lies in $\rH^{>0}(L)$, which is a proper ideal in $\rH^*(L)$, so it cannot be invertible.  We obtain the desired contradiction, proving the theorem.
\end{proof}

\begin{acknowledgements}
I am grateful to Noah Porcelli for useful correspondence and an anonymous referee for helpful comments.
\end{acknowledgements}

\bibliographystyle{siam} 
\bibliography{HLMbiblio} 

\end{document}